\documentclass[10pt]{article}
\usepackage{amsmath}
\usepackage{amsfonts}
\usepackage{amssymb}
\usepackage{amsthm}
\usepackage{geometry}                
\usepackage{graphicx}
\usepackage{epstopdf}
\usepackage{enumerate}
\begin{document}

\title{{\bf Periodic points of algebraic functions and Deuring's class number formula} }        
\author{Patrick Morton}        
\date{}          
\maketitle

\begin{abstract}  The exact set of periodic points in $\overline{\mathbb{Q}}$ of the algebraic function $\hat F(z)=(-1\pm \sqrt{1-z^4})/z^2$ is shown to consist of the coordinates of certain solutions $(x,y)=(\pi, \xi)$ of the Fermat equation $x^4+y^4=1$ in ring class fields $\Omega_f$ over imaginary quadratic fields $K=\mathbb{Q}(\sqrt{-d})$ of odd conductor $f$, where $-d =d_K f^2 \equiv 1$ (mod $8$).  This is shown to result from the fact that the $2$-adic function $F(z)=(-1+ \sqrt{1-z^4})/z^2$ is a lift of the Frobenius automorphism on the coordinates $\pi$ for which $|\pi |_2<1$, for any $d \equiv 7$ (mod $8$), when considered as elements of the maximal unramified extension $\textsf{K}_2$ of the $2$-adic field $\mathbb{Q}_2$.  This gives an interpretation of the case $p=2$ of a class number formula of Deuring.  An algebraic method of computing these periodic points and the corresponding class equations $H_{-d}(x)$ is given that is applicable for small periods.  The pre-periodic points of $\hat F(z)$ in $\overline{\mathbb{Q}}$ are also determined.
\end{abstract}

\section{Introduction.} 

In the papers \cite{d1} and \cite{d2} Deuring noted the following class number formulas. \bigskip

\noindent {\bf Deuring's Class Number Formulas.} \smallskip
$$\sum_{d_{p^f}}{h(d_{p^f})}+h_p = p^f, \ \ \textrm{for even} \ f=2,4,6 \cdots,$$
$$\sum_{d_{p^f}}{h(d_{p^f})}+2t_p-h_p=p^f, \ \ \textrm{for odd} \ f = 1, 3, 5, \cdots;$$
{\it where $h(d)$ is the class number of primitive quadratic forms of discriminant $-d$; $d_{p^f}$ runs over all positive integers, for which the principal form of discriminant $-d_{p^f}$ properly represents $p^f$; $h_p$ is the class number of the quaternion algebra $D_p=\mathbb{Q}_{\infty,p}$ which is ramified only at $p$ and the infinite prime $p_\infty$; and $t_p$ is the type number of $D_p$}. \bigskip

In these formulas, $h_p$ is the total number of $j$-invariants of supersingular elliptic curves in characteristic $p$, and $2t_p-h_p$ is the number of supersingular $j$-invariants which lie in the prime field $\mathbb{F}_p$ (see \cite{d2} and \cite{brm}, p. 97).  When these two numbers are the same, i.e. when all supersingular $j$-invariants lie in $\mathbb{F}_p$, then these formulas may be combined and ``inverted'' to give that
\begin{equation}
\sum_{d_{p^f}}{' \ h(d_{p^f})} = \sum_{k \mid f}{\mu(f/k)(p^k-h_p)} = \sum_{k \mid f}{\mu(f/k)p^k}, \ \ \textrm{for} \ f>1,
\end{equation}
where the primed sum is now taken over the positive integers, for which the principal form of discriminant $-d_{p^f}$ properly represents $p^f$ and {\it no smaller power of} $p$. \medskip

In \cite{m1} the formula (1) was reinterpreted (replacing $f$ with $n$) for the prime $p=3$ in the following form: Let $\mathfrak{D}_n^{(3)}$ be the set of discriminants $-d \equiv 1$ (mod $3$) of orders in imaginary quadratic fields $K=\mathbb{Q}(\sqrt{-d})$, for which the Frobenius automorphism $\displaystyle \tau_d=\left(\frac{\Omega_f/K}{\wp_3}\right)$, for a prime divisor $\wp_3$ of $3$ in $K$, has order $n$ in the Galois group $G=G(\Omega_f/K)$ of the ring class field $\Omega_f$ of conductor $f$ over $K$, where $-d=d_Kf^2$ and $d_K$ is the discriminant of $K$.  If $h(-d)$ is the class number of the order $\textsf{R}_{-d}$ of discriminant $-d$ in $K$, then \medskip
$$\sum_{-d \in \mathfrak{D}_n^{(3)}}{h(-d)} = \sum_{k \mid n}{\mu(n/k)3^k}, \ \ \textrm{for} \ n>1.$$

\noindent An independent proof of this formula was given, by interpreting the sum on either side as the number of periodic points of least period $n$ of a specific $3$-adic algebraic function defined and single-valued in a certain domain of the field $\textsf{K}_3$, the maximal unramified algebraic extension of the $3$-adic field $\mathbb{Q}_3$.  From the Artin reciprocity law we know that $\mathfrak{D}_n^{(3)}$ is the set of negative discriminants $-d \equiv 1$ (mod $3$) for which a prime ideal divisor $\wp_3$ of $3$ in the ring of integers $R_K$ of $K$ has order $n$ in the ring class group (mod $f$). \medskip

A similar interpretation for the prime $p=2$ was given in \cite{m2}, except that in formula (1) the prime $p$ was replaced by $2^2$:
\begin{equation}
\sum_{-d \in \mathfrak{D}'_n}{h(-d)} = \sum_{k \mid n}{\mu(n/k)2^{2k}}, \ \ \textrm{for} \ n>1;
\end{equation}
this is equivalent to the first half of the Deuring class number formula for the prime $p=2$.  Here $\mathfrak{D}'_n$ is the set of discriminants $-d \equiv 1$ (mod $8$), for which {\it the square} of the Frobenius automorphism $\displaystyle \tau_d=\left(\frac{\Omega_f/K}{\wp_2}\right)$ has order $n$ in the Galois group $G(\Omega_f/K)$ over $K=\mathbb{Q}(\sqrt{-d})$, where $-d=d_Kf^2$.  Once again, the number on either side of (2) was interpreted as the number of periodic points of least period $n$ of a specific $2$-adic algebraic function in a certain domain of $\textsf{K}_2$, the maximal unramified algebraic extension of $\mathbb{Q}_2$.  \medskip

In this note I show how the full formula (1) may be interpreted for the prime $p=2$ and $n>1$.  This arises from the fact that for ring class fields of a specific family of imaginary quadratic fields, the Frobenius automorphism $\tau_d$ can be represented by a {\it single} power series, independent of $d$, evaluated at one of a family of related generators for the fields $\Omega_f$.  \medskip

Before stating the precise result we recall the following definitions and results from \cite{lm}(Section 10) and \cite{m2}.
The Schl\"afli functions $\mathfrak{f}(\tau), \mathfrak{f}_1(\tau)$, $\mathfrak{f}_2(\tau)$ (see \cite{sch}, p. 148, or \cite{co1}, p. 256) can be defined by the infinite products
\begin{equation*}
\mathfrak{f}(\tau)=q^{-\frac{1}{48}} \prod_{n=1}^\infty{(1+q^{n-\frac{1}{2}})},\quad \mathfrak{f}_1(\tau)=q^{-\frac{1}{48}} \prod_{n=1}^\infty{(1-q^{n-\frac{1}{2}})},
\end{equation*}

\begin{equation*}
\mathfrak{f}_2(\tau)=\sqrt{2} \hspace{.05 in} q^{\frac{1}{24}} \prod_{n=1}^\infty{(1+q^n)}, \quad q = e^{2 \pi i \tau},
\end{equation*}
for $\tau$ in the upper half-plane $\mathbb{H}$.  Let $K=\mathbb{Q}(\sqrt{-d})$ be an imaginary quadratic field, for which $-d \equiv 1$ (mod $8$) is the discriminant of the order $\textsf{R}_{-d}$ in $K$, with conductor $f$, satisfying $-d=d_K f^2$.  Further, let $w \in K$ be given by
\begin{equation}
w=\frac{v+\sqrt{-d}}{2}, \quad v^2 \equiv -d \ (\textrm{mod} \ 16), \ \ v=1 \ \textrm{or} \ 3,
\end{equation}
and set
\begin{equation*}
a \equiv
\begin{cases}
\frac{-3d+5}{\ 16} \ (\textrm{mod} \ 4), &\textrm{if} \ v=3 \ \textrm{and} \ d \equiv 7 \ (\textrm{mod} \ 16),\cr
	\frac{-d+31}{16} \ (\textrm{mod} \ 4), &\textrm{if} \ v=1 \ \textrm{and} \ d \equiv 15 \ (\textrm{mod} \ 16).
\end{cases}
\end{equation*}
Then the numbers
\begin{equation}
\pi_d = i^a\frac{\mathfrak{f}_2(w/2)^2}{\mathfrak{f}(w/2)^2}, \ \ \xi_d=\frac{\beta}{2}=i^{-v}\frac{\mathfrak{f}_1(w/2)^2}{\mathfrak{f}(w/2)^2}
\end{equation}
lie in the ring class field $\Omega_f$ of conductor $f$ over $K$, and satisfy
\begin{equation*}
\pi_d^4+\xi_d^4=1.
\end{equation*}
The quantities $\pi_d$ and $\xi_d$ are conjugate algebraic integers over $\mathbb{Q}$ and $\Omega_f=\mathbb{Q}(\pi_d)=\mathbb{Q}(\xi_d)$.  Furthermore, if $\wp_2=(2,w)$ is a prime ideal divisor of $2$ in $K$, then $(2)=2R_K=\wp_2 \wp_2'$ in $K$, and we have
\begin{equation*}
(\pi_d)=\pi_d R_{\Omega_f} = \wp_2 R_{\Omega_f}, \ \ (\xi_d)=\xi_d R_{\Omega_f} = \wp_2'R_{\Omega_f}, \ \textrm{in} \ \Omega_f,
\end{equation*}
where $R_L$ denotes the ring of algebraic integers in the field $L$.  In addition, we will need the fact that there is an automorphism $\psi \in \textrm{Gal}(\Omega_f/\mathbb{Q})$ of order $2$ which interchanges $\pi_d$ and $\xi_d$ and therefore also interchanges the ideals $\wp_2$ and $\wp_2'$.  If $\tau \in \textrm{Gal}(\Omega_f/K)$, then $\tau^{-1}\psi \tau$ is an automorphism of order $2$ which interchanges $\pi=\pi_d^\tau$ and $\xi=\xi_d^\tau$.  \smallskip

Let $b_d(x)$ be the minimal polynomial over $\mathbb{Q}$ of the numbers $\pi_d$ and $\xi_d$.  Then $b_d(x)$ is a normal polynomial over $\mathbb{Q}$ (meaning that one of its roots generates a normal extension of $\mathbb{Q}$) and
\begin{equation*}
\textrm{deg}(b_d(x))=2h(-d),
\end{equation*} 
where $h(-d)$ is the class number of the order $\textsf{R}_{-d}$, i.e. the number of elements of the ideal class group of $\textsf{R}_{-d}$.  Recall from \cite{lm} that half of the roots of $b_d(x)$ are generators of the ideal $\wp_2 R_{\Omega_f}$ and half are generators of $\wp_2'R_{\Omega_f}$.  \medskip

With these definitions the following theorem holds. \bigskip

\noindent {\bf Theorem 1.} {\it Let $K=\mathbb{Q}(\sqrt{-d})$, where $-d \equiv 1$ (mod $8$), and let $d_K$ denote the discriminant of $K/\mathbb{Q}$.  Set $-d=d_Kf^2$, $\displaystyle \tau_d=\left(\frac{\Omega_f/K}{\wp_2}\right)$, and let $\pi$ be any root of the polynomial $b_d(x)$ for which $(\pi)=\wp_2$ in the ring of integers $R_{\Omega_f}$ of the ring class field $\Omega_f=\mathbb{Q}(\pi_d)=\mathbb{Q}(\pi)$ over $K$.  Further, let $F(z)$ be the algebraic function
\begin{equation}
F(z)=\frac{-1+\sqrt{1-z^4}}{z^2}=\sum_{n=1}^\infty{(-1)^n{\frac{1}{2} \atopwithdelims ( ) n} z^{4n-2}},
\end{equation}
defined for $z$ in the disc $\textsf{D}=\{z: |z|_2 <1\}$, a subset of the maximal unramified extension $\textsf{K}_2$ of the $2$-adic field $\mathbb{Q}_2$.  Then for any such integer $d$,
$$\pi^{\tau_d} = F(\pi) \ \ \textrm{in} \ \mathsf{K}_2,$$
if $\Omega_f \rightarrow \left(\Omega_f\right)_{\mathfrak{p}}$ is embedded in $\textsf{K}_2$ by completing at a prime divisor $\mathfrak{p}$ of $\wp_2$.}
\bigskip

This is an improvement and simplification over what I was able to show in \cite{m2}, since there I was only able to represent the action of $\tau_d^2$ by a power series evaluated at a generator of $\Omega_f$.  Note that Theorem 1 is analogous to the action of the polynomial $P(z)=z^k$ on cyclotomic fields $\mathbb{Q}(\zeta_n)$, where $\zeta_n=e^{2 \pi i/n}$ and $(n,k)=1$, since $\zeta_n \rightarrow P(\zeta_n)$ represents an automorphism for this family of abelian fields.  \medskip

Theorem 1 leads to the following result, with a substantially simpler proof than the proof that was given for the corresponding theorem in \cite{m2}. \bigskip

\noindent {\bf Theorem 2.} {\it (a) The periodic points of the $2$-adic function $F(z)$ in the disc $\textsf{D}=\{z: |z|_2 <1\} \subset \textsf{K}_2$ are $z=0$ and the roots $\pi$ of the polynomials $b_d(x)$ which lie in $\textsf{D}$, as $d$ runs over all positive integers $d \equiv 7$ (mod $8$)}.  \medskip

{\it (b) The periodic points of the multivalued function
$$\hat F(z)=\frac{-1\pm \sqrt{1-z^4}}{z^2},$$
satisfying $g(\hat F(z),z)=0$, with $g(x,y)=y^2 x^2+2x+y^2$, are $0, -1$, and the roots of the polynomials $b_d(z)=0$, as $d$ ranges over all positive integers $d \equiv 7$ (mod $8$).  This statement holds in any of the fields $\overline{\mathbb{Q}}_2, \overline{\mathbb{Q}}, \mathbb{C}.$} \bigskip

As in \cite{m2}, a periodic point of the multivalued algebraic function $\hat F(z)$ is defined to be a value $a$ (in an algebraically closed field $k$) for which there exist $n \in \mathbb{N}$ and $a_1, \dots, a_{n-1} \in k$, for which the minimal polynomial $g(x,z)$ of $x=F(z)$ over $k(z)$ satisfies
$$g(a,a_1)=g(a_1,a_2)=\cdots =g(a_i,a_{i+1})=\dots=g(a_{n-1},a)=0.$$

Theorem 2 shows again, as in \cite{m2}, that all ring class fields of odd conductor $f$ over fields $K=\mathbb{Q}(\sqrt{-d})$ with $-d \equiv 1$ (mod $8$) can be generated over $\mathbb{Q}$ by individual periodic points of the algebraic function $\hat F(z)$; moreover, that all periodic points of $\hat F(z)$, with the exception of $z=0,-1$, generate ring class fields over fields $K$ in this same family. \bigskip

\noindent {\bf Corollary.}  {\it If $\mathfrak{D}_n=\mathfrak{D}_n^{(2)}$ is the set of negative discriminants $-d \equiv 1$ (mod $8$) for which the Frobenius automorphism $\tau=\left(\frac{\Omega_f/K}{\wp_2}\right)$ has order $n$ in $Gal(\Omega_f/K)$, with $K=\mathbb{Q}(\sqrt{-d})$, $-d=d_K f^2$ and $2 \cong \wp_2 \wp_2'$ in $R_K$, then for any $n > 1$ we have the class number formula}
\begin{equation}
\sum_{-d \in \mathfrak{D}_n}{h(-d)}=\sum_{k|n}{\mu(n/k)2^k}.
\end{equation} \medskip

This corollary is a consequence of Theorem 2(a) and the fact that the period $n$ of a periodic point $\pi \in \textsf{D}$ of $F(z)$ is the order of the automorphism $\tau_d$, by Theorem 1.  Thus, the sum in the corollary is the number of periodic points of $F(z)$ in $\textsf{D}$ with primitive (i.e., minimal) period $n$.  \medskip

In \cite{m2} the analogue of Theorem 2 was proved for the algebraic function
$$F_1(z)=-\frac{1+\sqrt[4]{1-z^4}}{1-\sqrt[4]{1-z^4}}=-\frac{(1+\sqrt[4]{1-z^4})^2(1+\sqrt{1-z^4})}{z^4}.$$
Thus, the discussion here proves that $F(z)=-z^{-2}+\left(\frac{F_1(z)+1}{zF_1(z)-z}\right)^2$ and $F_1(z)$ have the same periodic points.  \medskip

The function $F(z)$ is closely related to the function $T(z)$ which is defined as follows.  There is an isogeny
$$ \phi: E_\lambda \rightarrow E_{\lambda_1}$$
of degree $2$ from the Legendre normal form
$$E_\lambda: \ Y^2=X(X-1)(X-\lambda)$$
for the parameter $\lambda$ to the Legendre normal form for the parameter $\lambda_1$, and the formula for $\lambda_1=T(\lambda)$ is
$$\lambda_1=\frac{(1 - \sqrt{1-\lambda})^4}{\lambda^2}.$$
The function $T(z)$ defined by
$$T(z)=\frac{(1-\sqrt{1-z})^4}{z^2}=\frac{1}{z^2}\left(\sum_{n \ge 1}{{\frac{1}{2} \atopwithdelims ( ) n}(-1)^{n+1}z^n}\right)^4,$$
for $z$ in the disc $\textsf{D}=\{z: |z|_2 < 1\} \subset \textsf{K}_2$, is related to $F(z)$ by the formula
$$T(z^4) = F(z)^4.$$
This yields the following theorem for the periodic points of the multivalued function $\hat T(z)=\frac{(1 \pm \sqrt{1-z})^4}{z^2}$.

\bigskip

\noindent {\bf Theorem 3.} {\it The periodic points of the multivalued algebraic function $y=\hat T(z)=\frac{(1 \pm \sqrt{1-z})^4}{z^2}$ defined by the equation
$$C: \ \ \tilde g(z,y)=z^2 y^2-2(z^2-8z+8)y+z^2=0$$
are the numbers in the set
$$S=\{0,1\} \cup \{\xi^4: \exists d > 0, d\equiv 7 \ (\textrm{mod} \ 8) \ s.t. \ b_d(\xi)=0\}.$$
These are the numbers
$$\xi^4=\frac{\mathfrak{f}_1(w/2)^8}{\mathfrak{f}(w/2)^8}=1-\lambda\left(\frac{w}{2}\right), \ \ \pi^4=\frac{\mathfrak{f}_2(w/2)^8}{\mathfrak{f}(w/2)^8}=\lambda\left(\frac{w}{2}\right)$$
and their conjugates over $\mathbb{Q}$, where the number $w$ has the form}
$$w=\frac{v+\sqrt{-d}}{2}, \  v^2 \equiv -d \ (\textrm{mod} \ 16), \ d \equiv 7 \ (\textrm{mod} \ 8).$$

The function $\lambda(z)$ in this theorem is the classical $\lambda$-function, which is a modular function for the principal congruence group $\Gamma[2]$.  (See \cite{cha} and \cite{scho}.)  Theorems 2 and 3 give two examples of algebraic functions, whose periodic points are values of modular functions.  \medskip

In Section 3 I show how to use the simple polynomial $g(x,y)=x^2 y^2+2x+y^2$ and iterated resultants to compute the minimal polynomials $b_d(x)$ in a purely algebraic way.  I show in Sections 3 and 4 that the particular polynomials $b_d(x)$, for which $d \in \mathfrak{D}_k$ and $k \mid n$, together with $x$ and $x+1$, make up the exact set of irreducible factors of an $(n-1)$-fold iterated resultant defined using $g(x,y)$.  It seems quite remarkable that the minimal polynomials of values of modular functions can be found in this way.  In particular, this gives an algebraic method for computing generators of the ring class fields $\Omega_f$ of fields of the type $K=\mathbb{Q}(\sqrt{-d})$, with $-d \equiv 1$ (mod $8$), and therefore a purely algebraic method for computing the corresponding class equations $H_{-d}(x)$.  See Theorem 9 and Tables 1, 2, and 3 at the end of the paper. \medskip

The periodic points ($\neq 0,1$) of the function $\hat T$ in Theorem 3 generate ring class fields of odd conductor over the quadratic fields of the form $K=\mathbb{Q}(\sqrt{-d})$, $-d \equiv 1$ (mod $8$), as was proved in \cite{lm}.  In Section 4 I use an identity for the modular function $\lambda(z)$ to show that the pre-periodic points of $\hat T$ (of level $r \ge 2$, see Section 5) generate ring class fields of even conductor over fields $K$ in the same family; and conversely, every ring class field of even conductor over such a field $K$ is generated over $K$ by a pre-periodic point of $\hat T$.  This result, summarized in Theorem 13, proves Conjecture 2 in the paper \cite{m2} for the prime $p=2$.  The discussion here also gives an alternate proof of Conjecture 1 (for $p=2$) in that paper.  Similar results holds for the periodic and pre-periodic points of $\hat F$, as we show in Theorems 14 and 15.  In particular, the collection of ring class fields over fields $K$ in this family coincides with the collection of normal closures over $\mathbb{Q}$ of fields generated by individual periodic or pre-periodic points of the algebraic function $\hat F$.

\section{Lifting the Frobenius automorphism on roots of $b_d(x)$.}
In this section $\pi=\pi_d^\sigma$ will be any root of $b_d(x)=0$ which is conjugate to $\pi_d$ over $K=\mathbb{Q}(\sqrt{-d})$, and $\xi=\xi_d^\sigma$ will be the root of $b_d(x)=0$ for which $\pi^4+ \xi^4=1$.  See equations (3) and (4). Changing notation slightly, we let $\psi \in \textrm{Gal}(\Omega_f/\mathbb{Q})$ be the automorphism of order $2$ for which
$$\pi^\psi=\xi, \ \ \wp_2^\psi=\wp_2'.$$

We recall the following ideal factorizations from \cite{lm}(Theorem 8.6 and proof). \bigskip

\noindent {\bf Lemma 4.} {\it If $\cong$ denotes divisor equality; $2 \cong \wp_2 \wp_2'$ in the ring of integers $R_K$ of $K=\mathbb{Q}(\sqrt{-d})$; and $\beta=2\xi$; then $\beta \cong \wp_2 \wp_2'^2$ and
$$\beta-2 \cong \wp_2^2 \wp_2', \ \ \beta+2 \cong \wp_2^3 \wp_2', \ \ \beta^2+4 \cong \wp_2^3 \wp_2'^2$$
in the ring $R_{\Omega_f}$}.  \bigskip 

\noindent {\bf Lemma 5.} {\it Using the notation of Theorem 1, we have}
$$1+\pi^2 \cong \wp_2' \ \ \textrm{and} \ \ 1-\pi^2 \cong \wp_2'^3,$$
{\it so that}
$$\frac{(1-\pi^2)^2}{(1+\pi^2)^2} \cong \wp_2'^4.$$
\smallskip

\noindent {\it Proof.} From \cite{lm} we have
$$(1+\pi^2)(1-\pi^2)=1-\pi^4 = \xi^4 \cong \wp_2'^4.$$
Furthermore, using Lemma 4, we have
$$1-\xi^2 = (1-\xi)(1+\xi) = \frac{(2-\beta)(2+\beta)}{4} \cong \frac{\wp_2^2 \wp_2' \wp_2^3 \wp_2'}{\wp_2^2 \wp_2'^2} = \wp_2^3.$$
Now apply the automorphism $\psi \in G(\Omega_f/K)$ which switches the numbers $\pi$ and $\xi$ and the ideals $\wp_2$ and $\wp_2'$: this gives $1-\pi^2 = (1-\xi^2)^\psi \cong \wp_2'^3$ and verifies the assertions. $\square$ \bigskip

\noindent {\bf Theorem 6.} {\it If $\tau = \left(\frac{\Omega_f/K}{\wp_2}\right)$, we have}
$$\xi^{4\tau^{-1}}=(\xi^{\tau^{-1}})^4 = \frac{(1-\pi^2)^2}{(1+\pi^2)^2}.$$
\smallskip

\noindent {\it Proof.} Letting $\alpha$ denote a solution of $16\alpha^4+16\beta^4=\alpha^4\beta^4$, we have as in \cite{lm} (pp. 1967-68), for a suitable basis quotient $w$ of an ideal $\mathfrak{a}$ (prime to $f$), that
$$j\left(\frac{w}{2}\right)=\frac{(\alpha^8-16\alpha^4+256)^3}{\alpha^8(\alpha^4-16)^2}=\frac{(\beta^8-16\beta^4+256)^3}{\beta^8(\beta^4-16)^2}=\frac{256(\xi^8-\xi^4+1)^3}{\xi^8(\xi^4-1)^2}=J(\xi^4),$$
where
$$J(x)=\frac{256(x^2-x+1)^3}{x^2(x-1)^2}.$$
Since $j(w/2)^{\tau^{-1}}=j(w/4)$, we find that
$$j\left(\frac{w}{4}\right)=j\left(\frac{w}{2}\right)^{\tau^{-1}}=J(\xi^{4\tau^{-1}}).$$
On the other hand, straightforward calculation shows that
$$J\left(\frac{(1-\pi^2)^2}{(1+\pi^2)^2}\right)=\frac{16(\pi^8+14\pi^4+1)^3}{\pi^4(\pi^4-1)^4}.$$
Replacing $\pi^4$ by $1-\xi^4$ in the last expression yields
$$J\left(\frac{(1-\pi^2)^2}{(1+\pi^2)^2}\right)=\frac{16(\xi^8-16\xi^4+16)^3}{\xi^{16}(1-\xi^4)}.$$
Using $\xi=\beta/2$, this and \cite{lm}(eq. (6.2)) yield
$$J\left(\frac{(1-\pi^2)^2}{(1+\pi^2)^2}\right)=\frac{(\beta^8-256\beta+4096)^3}{\beta^{16}(16-\beta^4)}=j\left(\frac{w}{4}\right),$$
which gives that
$$J(\xi^{4\tau^{-1}})=J\left(\frac{(1-\pi^2)^2}{(1+\pi^2)^2}\right).$$
Setting $z_1=\xi^{4\tau^{-1}}$ and $z_2=\frac{(1-\pi^2)^2}{(1+\pi^2)^2}$, this implies that $z_1$ and $z_2$ are related by an element of the anharmonic group:
$$z_2 \in \{\frac{1}{z_1}, \frac{z_1}{z_1-1}, \frac{z_1-1}{z_1}, \frac{1}{1-z_1}, 1-z_1, z_1\}.$$
Now, $z_2$ cannot be equal to any of the first four elements, since $z_1\cong z_2 \cong \wp_2'^4$ and $1-z_1 \cong \wp_2^4$ imply that these four elements are not integral.  Similarly, $z_2 \neq 1-z_1$, forcing $z_2=z_1$.  This proves the theorem.  $\square$  \bigskip

Applying the automorphism $\tau^{-1}$ to the equation $\pi^4+\xi^4=1$ gives
$$\pi^{4\tau^{-1}}+\frac{(1-\pi^2)^2}{(1+\pi^2)^2}=1;$$
hence, the points $(x,y)=(\pi,\pi^{\tau^{-1}})$ and $(x,y) = (\pi^\tau,\pi)$ satisfy the equation
\begin{equation}
y^4+\frac{(1-x^2)^2}{(1+x^2)^2}-1=\frac{(x^2 y^2+2x+y^2)(x^2 y^2-2x+y^2)}{(x^2+1)^2}=0.
\end{equation}
Note that the first factor in the numerator is the polynomial $g(x,y)=x^2 y^2+2x+y^2$ defined in Theorem 2 above, and that the function $F(z)$ defined in (5) satisfies $g(F(z),z)=0$.  \bigskip

We expand $F(z)=\frac{-1+\sqrt{1-z^4}}{z^2}$ in a $2$-adic series in the field $\textsf{K}_2$:
$$F(z)=\sum_{n=1}^\infty{(-1)^n{\frac{1}{2} \atopwithdelims ( ) n} z^{4n-2}}.$$
If $|z|_2 < 1$, then $z=2z_1$, with $|z_1|_2 \le 1$, so that
$$F(z)=F(2z_1)=\sum_{n=1}^\infty{(-1)^n 2^{4n-2}{\frac{1}{2} \atopwithdelims ( ) n} z_1^{4n-2}}$$
$$ \hspace{.9 in} = -\sum_{n=1}^\infty{2^{2n-1}C_{n-1} z_1^{4n-2}}=-2z_1^2-8z_1^6-\cdots,$$
where $C_{n-1}=(-1)^{n-1} 2^{2n-1}{\frac{1}{2} \atopwithdelims ( ) n} \in \mathbb{Z}$ is the Catalan number.  Hence, the series converges for $z$ in the disc $\textsf{D}=\{z: |z|_2<1\}$, and maps this disc into itself.  This allows us to iterate the function $F(z)$ on $\textsf{D}$.  For the proof of Theorem 1 we also need the relation
\begin{equation}
\pi^{\tau^2}=\frac{\xi+1}{\xi-1}
\end{equation}
from \cite{lm}(Prop. 8.5). \medskip

\noindent {\it Proof of Theorem 1.} \medskip

We would like to see that $\pi^\tau=F(\pi)$.  From (7) we know that the point $(x,y)=(\pi^\tau,\pi)$ lies either on the curve $x^2 y^2+2x+y^2=0$ or on $x^2 y^2-2x+y^2=0$.  Suppose this point lies on the latter curve; then
$$\pi^{2\tau} \pi^2-2\pi^\tau+\pi^2=0$$
implies that
$$\pi^\tau=\frac{1-\sqrt{1-\pi^4}}{\pi^2}=\frac{1 - \xi^2}{\pi^2},$$
since the expression $\frac{1+\xi^2}{\pi^2}\cong \wp_2^{-1}$ is not integral.  This makes use of the calculations in Lemma 4, according to which $1+\xi^2 \cong \wp_2 \cong \pi$.  Applying $\tau$ to the last displayed equation and using (8) gives
$$\frac{\xi+1}{\xi-1}=\pi^{\tau^2}=\frac{1-\xi^{2\tau}}{\pi^{2\tau}}=\frac{1-\xi^{2\tau}}{(1-\xi^2)^2/\pi^4},$$
so that
$$\frac{(\xi+1)^3(\xi-1)}{\pi^4}=1-\xi^{2\tau}.$$
Using $\pi^4=1-\xi^4$ in the last relation yields
$$1-\xi^{2\tau}=\frac{(\xi+1)^3(\xi-1)}{1-\xi^4}=-\frac{(\xi+1)^2}{1+\xi^2}$$
and
$$\xi^{2\tau}=1+\frac{(\xi+1)^2}{1+\xi^2}=\frac{2(1+\xi+\xi^2)}{1+\xi^2}.$$
Since $1+\xi^2 \cong \wp_2$ and $\xi \cong \wp_2'$, it follows that $(1+\xi+\xi^2,2)=1$.  Thus, the right hand side in the last displayed equation is $\mathfrak{a} \wp_2'$, where $(\mathfrak{a},2)=1$.  However, the left side is $\xi^{2\tau} \cong \wp_2'^2$, giving a contradiction.  This proves that
\begin{equation}
\pi^\tau = \frac{-1+\sqrt{1-\pi^4}}{\pi^2}=\frac{-1+\xi^2}{\pi^2}=F(\pi).
\end{equation}
$\square$ \medskip

This proves Theorem 1.  Iterating equation (9), and noting that $\tau$ can be viewed as an automorphism of the local extension $\mathbb{Q}_2(\pi)/\mathbb{Q}_2$, we find that
$$\pi^{\tau^n}=F^n(\pi), \ \ \textrm{for} \ n \ge 1.$$
Therefore, we have the following result. \bigskip

\noindent {\bf Theorem 7.} {\it The roots $\pi$ of the polynomials $b_d(x)$, for which $\pi \cong \pi_d \cong \wp_2$, are all periodic points in $\textsf{D} \subset \textsf{K}_2$ of the algebraic function $F(z)$.  The period $n$ of such a number $\pi$ is the order of $\tau_d$ in the Galois group of the ring class field $\Omega_f$ over $\mathbb{Q}$.}\bigskip

Since half of the roots of $b_d(x)=0$ are conjugate to $\pi_d \cong \wp_2$ over $K=\mathbb{Q}(\sqrt{-d})$ and half are conjugate to $\xi_d \cong \wp_2'$, this shows that there are $h(-d)$ periodic points of $F(z)$ in the disk $\textsf{D}$ coming from the roots of $b_d(x)=0$, for a given $d$.  \medskip

Applying the automorphism $\psi$ to (9) and using $\tau \psi=\psi \tau^{-1}$ implies
\begin{equation}
\xi^{\tau^{-1}} = \frac{-1+\pi^2}{\xi^2}.
\end{equation}
This agrees with the result of Theorem 1, since
$$\frac{(-1+\pi^2)^4}{\xi^8}=\frac{(\pi^2-1)^4}{(\pi^4-1)^2}=\frac{(1-\pi^2)^2}{(1+\pi^2)^2}.$$
Moreover, $(\xi,\xi^\tau)$ is also a point on the curve $x^2 y^2+2x+y^2=0$, as can be seen by applying the automorphism $\psi \tau$ to the point $(\pi^\tau,\pi)$:
$$(\pi^\tau,\pi)^{\psi \tau}=(\pi^{\tau \psi \tau},\pi^{\psi \tau})=(\pi^{\psi \tau^{-1} \tau},\xi^\tau)=(\xi,\xi^\tau).$$
It follows that $\xi$ is a periodic point of the inverse algebraic function $F^{-1}(z)$, for which $g(z,F^{-1}(z))=0$.  Hence, the remaining $h(-d)$ roots $\xi$ of $b_d(x)=0$ are all periodic points of $F^{-1}(z)$.  See equation (14) below.  \medskip

\section{Iterated resultants.}

We turn now to the proof that the roots of the polynomials $b_d(x)$, together with $0$ and $-1$, are the {\it only} periodic points of the multivalued function $\hat F(z)$.  We let
$$g(x,y)=x^2 y^2+2x+y^2$$
as before, and define
\begin{align*}
R^{(1)}(x,x_1)&=g(x,x_1),\\
R^{(2)}(x,x_2)&=\textrm{Res}_{x_1}(g(x,x_1),g(x_1,x_2));
\end{align*}
and inductively,
$$R^{(n)}(x,x_n)=\textrm{Res}_{x_{n-1}}(R^{(n-1)}(x,x_{n-1}),g(x_{n-1},x_n)), \ n \ge 2.$$
Putting $x_n=x$ gives the polynomial
$$R_n(x)=R^{(n)}(x,x), \ \ n \ge 1.$$
The roots of $R_n(x)$ are exactly the elements $a \in \mathbb{A}_2=\overline{\mathbb{Q}}_2$ for which there are $a_1, \dots, a_{n-1} \in \mathbb{A}_2$ satisfying the simultaneous equations
\begin{equation}
g(a,a_1)=g(a_2,a_3)=\cdots =g(a_{n-1},a)=0;
\end{equation}
i.e. the $a$'s are exactly the periodic points of $\hat F(z)$ with period $n$.  \medskip

The same arguments as in \cite{m2} lead to the factorizations
$$R_n(x)=\prod_{k|n}{\textsf{P}_k(x)},$$
\begin{equation}
\textsf{P}_n(x)=\prod_{k|n}{R_k(x)^{\mu(n/k)}}.
\end{equation}
where $\mu(n)$ is the M\"obius $\mu$-function.  This is done by relating the polynomials $R_n(x)$ and $\textsf{P}_n(x)$ to the corresponding polynomials 
$\tilde R_n(x)$ and $\widetilde{\textsf{P}}_n(x)$ obtained by replacing $g(x,y)$ in the above definitions by the polynomial
$$g_1(x,y)=\frac{g(2x,2y)}{4}=4x^2 y^2+x+y^2, \ \ g(x,y)=x^2 y^2+2x+y^2.$$
Note that
$$g_1(x,y) \equiv y^2+x \ (\textrm{mod} \ 2).$$
It follows easily by induction that
$$\tilde R^{(n)}(x,x_n) \equiv x_n^{2^n}+x \ (\textrm{mod} \ 2), \ \ n \ge 1,$$
and therefore
$$\tilde R_n(x) \equiv x^{2^n}+x \ (\textrm{mod} \ 2), \ \ n \ge 1.$$
Hensel's Lemma implies that $\tilde R_n(x)$ has at least $2^n$ distinct roots in $\textsf{K}_2$, of which $2^n-1$ are units in $\textsf{K}_2$.  It can also be checked that
$$R_n(2x)=2^{2^n} \tilde R_n(x),$$
which implies that $R_n(x)$ also has at least $2^n$ distinct roots in $\textsf{K}_2$, of which $2^n-1$ are prime elements in $\textsf{K}_2$.  Note that $x=0$ is certainly a root of $R_n(x)$ for any $n$.  \medskip

It follows from the identity
\begin{equation}
(x-1)^2(y-1)^2g\left(\frac{x+1}{x-1},\frac{y+1}{y-1}\right)=4g(y,x)
\end{equation}
that for every root $a \in \textsf{K}_2$ of $R_n(x)$ the quantity $b=\frac{a+1}{a-1} \in \textsf{K}_2$ is also a root.  This is because $g(a,a_1)=g(a_1,a_2)=\cdots =g(a_{n-1},a)=0$ and $b_i=\frac{a_i+1}{a_i-1}$ imply that
$$g(b,b_{n-1})=\cdots=g(b_2,b_1)=g(b_1,b)=0.$$
The roots $b$ are distinct from all the roots $a$, since the $b$'s are all units.  Hence, $R_n(x)$ has $2^{n+1}$ distinct roots in $\textsf{K}_2$.  It is not hard to see that $\textrm{deg}(R_n(x))=2^{n+1}$, so this accounts for all the roots.  (See the proof of the Lemma in \cite{m2}, pp.727-728.) Furthermore, if $k \mid n$, then roots of $R_k(x)$ are also roots of $R_n(x)$.  It follows that the expression $\textsf{P}_n(x)$ defined in (12) is a polynomial.  This gives that
$$\textrm{deg}(R_n(x)) = 2^{n+1}, \ \ \textrm{deg}(\textsf{P}_n(x))=2\sum_{k|n}{\mu(n/k)2^k}.$$
The roots of the polynomial $\textsf{P}_n(x)$ are exactly the periodic points of $\hat F(z)$ of minimal period $n$. \medskip

This discussion proves: \bigskip

\noindent {\bf Theorem 8.} {\it All the periodic points of the multivalued algebraic function $\hat F$ lie in the maximal unramified, algebraic extension $\textsf{K}_2$ of the $2$-adic field $\mathbb{Q}_2$.} \bigskip

Irreducible factors of the polynomials $\textsf{P}_n(x)$ are listed in Tables 1, 2, and 3 for small values of $n$. \medskip

The identity (13) also implies
$$g\left(\frac{z+1}{z-1},\frac{F(z)+1}{F(z)-1}\right)=0,$$
so that
\begin{equation}
F^{-1}\left(\frac{z+1}{z-1}\right)=\frac{F(z)+1}{F(z)-1},
\end{equation}
where $F^{-1}(z)$ is defined by
$$F^{-1}(z)=-\frac{\sqrt{-2z(z^2+1)}}{z^2+1}.$$
Equation (14) shows that $F^{-1}(z)$ can be defined as a single valued function on the image of the disc $\textsf{D}$ under the map $\phi(z)=\frac{z+1}{z-1}$.  This image is the set
$$\phi(\textsf{D})=\{z  :  |z+1|_2 \le 2^{-2} \} =\{z : z=4w-1, w \in \mathfrak{o}_2 \},$$
where $\mathfrak{o}_2$ is ring of integers in $\textsf{K}_2$.
\medskip

From the preceding discussion we also see that
\begin{equation}
\widetilde{\textsf{P}}_n(x) \equiv \prod_{k|n}{(x^{2^k}+x)^{\mu(n/k)}} \ (\textrm{mod} \ 2),
\end{equation}
where the right side is the product of the irreducible polynomials of degree $n$ in $\mathbb{F}_2[x]$.  This implies that over $\mathbb{Q}_2$, the irreducible factors of $\widetilde{\textsf{P}}_n(x)$, and hence also of $\textsf{P}_n(x)$, have degree $n$.  If $a$ is a periodic point of $\hat F$ and $a_1, \dots, a_{n-1}$ are the associated elements of $\textsf{K}_2$ satisfying (11), then the $a_i$ are also roots of $R_n(x)$, as can be seen by cyclically permuting the equations in (11).  Hence, the roots of $R_n(x)$ consist of complete orbits under $\hat F$.  The same holds for the polynomial $\tilde R_n(x)=2^{-2^n}R_n(2x)$ under the conjugate map $\tilde F(z)=\frac{F(2z)}{2}$.  Since
\begin{equation}
\frac{F(2z)}{2} = -z^2-4z^6 + \cdots \equiv z^2 \ (\textrm{mod} \ 2), \ \ \textrm{for} \ |z|_2 \le 1,
\end{equation}
and the Frobenius map $z \rightarrow z^2$ fixes the irreducible factors of degree $n$ over $\mathbb{F}_2$, it follows that the roots of an irreducible factor of $\textsf{P}_n(x)$ over $\mathbb{Q}_2$ consist of: a complete orbit under the map $F(z)$, if those roots lie in $\textsf{D}$; and a complete orbit under $F^{-1}(z)$, if those roots lie in $\phi(\textsf{D})$.
\medskip

Finally, note that Theorem 7 implies that $b_d(x) \mid \textsf{P}_n(x)$ whenever $\textrm{ord}(\tau_d)=n$.  In the next section we show that the polynomials $b_d(x)$, together with $x$ and $x+1$, are the only irreducible factors of $\textsf{P}_n(x)$.  See Theorem 9 below.  \medskip

\section{Elliptic curves and periodic points.}

Consider the isogeny $\phi_1$ of degree $2$ on the elliptic curve
$$E_\lambda: \ y^2=x(x-1)(x-\lambda)$$
which is induced by the translation map
$$(x,y)^\rho=(x,y)-(1,0)=\left(\frac{x-\lambda}{x-1},\frac{\lambda-1}{(x-1)^2}y\right);$$
namely
$$\phi_1(x,y)=(u,v),$$
where
\begin{align*}
u &= x+x^\rho=x+\frac{x-\lambda}{x-1}=\frac{x^2-\lambda}{x-1},\\
v &= y+y^\rho =y+\frac{\lambda-1}{(x-1)^2}y=\frac{x^2-2x+\lambda}{(x-1)^2}y.
\end{align*}
The image of $E_\lambda$ under this isogeny is the curve
$$E_1: \ v^2=(u-\lambda)(u^2-4u+4\lambda).$$
Replacing $u$ by $u_1+\lambda$ gives the equation
$$v^2=u_1(u_1^2+(2\lambda-4)u_1+\lambda^2)=u_1(u_1-\gamma)(u_1-\gamma'),$$
where the roots of the quadratic are
$$\gamma, \gamma'=-\lambda+2 \pm 2\sqrt{1-\lambda}=(1 \pm \sqrt{1-\lambda})^2.$$
Now setting $X=\frac{u_1}{\gamma}$ and $Y=\frac{v}{\gamma^{3/2}}$ yields the curve
$$E_{\lambda_1}: \ Y^2=X(X-1)(X-\frac{\gamma'}{\gamma}),$$
where
$$\lambda_1=\frac{\gamma'}{\gamma}=\frac{(1 - \sqrt{1-\lambda})^2}{(1 + \sqrt{1-\lambda})^2}=\frac{(1- \sqrt{1-\lambda})^4}{\lambda^2}.$$
\medskip

These transformations yield an isogeny $\phi: E_\lambda \rightarrow E_{\lambda_1}$, with
$$\lambda_1=\frac{\gamma'}{\gamma}=\frac{(1- \sqrt{1-\lambda})^4}{\lambda^2}=T(\lambda),$$
for which $\phi(0,0) = (0,0)$.  The functions $T(z)=\frac{(1- \sqrt{1-z})^4}{z^2}$ and $F(z)$ satisfy the relationship
\begin{equation}
T(z^4)=\frac{(1- \sqrt{1-z^4})^4}{z^8}=F(z)^4,
\end{equation}
which implies that
$$T^n(z^4)=F^n(z)^4, \ \ n \ge 1, \ \ |z|_2 < 1.$$
If $\pi$ is any periodic point of $F(z)$ in the disk $\textsf{D}$, with period $n$, then
$$T^n(\pi^4)=F^n(\pi)^4=\pi^4$$
shows that $\pi^4$ is a periodic point of $T$ of period $n$ in $\textsf{D}$.  Conversely, if $\pi^4 \in \textsf{D}$ is a periodic point of $T(z)$, then $\pi \in \textsf{D}$ and $T^n(\pi^4)=\pi^4=F^n(\pi)^4$ implies that $F^n(\pi)=\pm \pi$, since $\sqrt{-1} \notin \textsf{K}_2$.  If $F^n(\pi)=-\pi$, then $F(-z)=F(z)$ implies that $-\pi$ is a periodic point of $F(z)$.  If $\pi$ were also a periodic point of $F(z)$, with period $m$, then $F^m(\pi)=\pi$ implies
$$\pi = F^{\circ mn}(\pi)=(F^{\circ n})^{\circ m}(-\pi) = -\pi,$$
giving $\pi=0$.  (The composition symbol $\circ$ is included here to emphasize that the powers are compositions.) Thus, periodic points $\pi^4$ of $T(z)$ in $\textsf{D}$ are in 1-1 correspondence with the periodic points $\pi$ of $F(z)$ in $\textsf{D}$.
\bigskip

{\it Proof of Theorem 2.}
Let $a$ be any periodic point of $\hat F(z)$ in $\textsf{K}_2$ with primitive period $n>1$; then the discussion in Section 3 shows that $a \in \textsf{D} \cup \phi(\textsf{D})$.  By replacing $a$ by $\phi(a)$ we assume $a \in \textsf{D}$.  Then $g(a_1,a)=0$ implies that $a_1 \in \textsf{D}$, so that $a_1=F(a)$.  Let $F^i(a)=a_i$ for $1 \le i \le n-1$, and $F^n(a)=a=a_0=a_n$ with $n$ smallest.  Each of the quantities $a_i$ is a periodic point of $F(z)$, so each $a_i^4$ is a periodic point of $T(z)$ in $\textsf{D}$.  Moreover,
$$T(a_i^4)=a_{i+1}^4, \ 0 \le i \le n-1.$$
Thus, for each $i$ with $0 \le i \le n-1$ there is an isogeny
$$\phi_i:E_{a_i^4} \rightarrow E_{a_{i+1}^4}.$$
Hence, $\iota=\phi_{n-1} \circ \phi_{n-2} \circ \cdots \circ \phi_1 \circ \phi_0$ is an isogeny from $E_{a^4}$ to itself.  Moreover, by the above discussion, each $\phi_i$ takes $(0,0)$ to $(0,0)$, so that $\iota$ fixes $(0,0) \in E_{a^4}$.  This implies that $\iota$ is a cyclic isogeny of degree $2^n$, and therefore the $j$-invariant of $E_{a^4}$, namely
$$j(E_{a^4})=2^8\frac{(a^8-a^4+1)^3}{a^8(a^4-1)^2},$$
is a root of the modular equation
$$\Phi_{2^n}(x,x)=c_n \prod_{-d}{H_{-d}(x)^{r(d,2^n)}},$$
where the product is over orders $\textsf{R}_{-d}$ of discriminant $-d$ in imaginary quadratic fields, $H_{-d}(x)$ is the class polynomial of discriminant $-d$, and
$$r(d,m)=|\{\lambda \in \textsf{R}_{-d}: \lambda \ \textrm{primitive}, N(\lambda)=m\}/\textsf{R}_{-d}^\times|.$$
See \cite{co1}.  Therefore, $x=a^4$ is a root of the polynomial
$$L_d(x)=(x^2-x)^{2h(-d)}H_{-d}\left(\frac{2^8(x^2-x+1)^3}{x^2(x-1)^2}\right).$$
Now the argument of \cite{m2}, pp. 736-737, applies word for word (with $\xi$ replaced by $a$), and shows that $a$ is a root of a polynomial $b_d(x)$, where $d \equiv 7$ (mod $8$) and $-d \in \mathfrak{D}_n$.  This proves Theorem 2(a).  Theorem 2(b) is immediate from the fact that the irreducible factors of the polynomial $R_n(x)$ are the $b_d(x)$, independent of which field $\overline{\mathbb{Q}}_2, \overline{\mathbb{Q}}, \mathbb{C}$ we are working in.  $\square$  \bigskip

The above arguments imply the following result concerning the polynomial $\textsf{P}_n(x)$ in (12). \bigskip

\noindent {\bf Theorem 9.} {\it The polynomial $\textsf{P}_n(x)$ is given by the formula
$$\textsf{P}_n(x)=\prod_{-d \in \mathfrak{D}_n}{b_d(x)}, \ \ \textrm{for} \ n>1,$$
where $\mathfrak{D}_n$ is defined (as in the Corollary to Theorem 2) as the set of negative discriminants $-d \equiv 1$ (mod $8$) for which $\tau_d=\left(\frac{\Omega_f/K}{\wp_2}\right)$ has order $n$ in $\textrm{Gal}(\Omega_f/K)$.  In particular, equating degrees yields}
\begin{equation}
2 \sum_{k \mid n}{\mu(n/k)2^k}=2 \sum_{-d \in \mathfrak{D}_n}{h(-d)}, \ \  \textrm{for} \ n>1 .
\end{equation}

For $n=1$ we have
\begin{equation}
\textsf{P}_1(x)=R_1(x)=g(x,x)=x(x+1)(x^2-x+2)=x(x+1)b_7(x).
\end{equation}

Theorem 9 shows that the polynomials $b_d(x)$ may be computed as the irreducible factors of the iterated resultants $R_n(x)$ in Section 3.  Factoring $R_n(x)$ for a fixed $n$ yields the complete set of polynomials $b_d(x)$, for $d \in \mathfrak{D}_n$.  See Tables 1, 2, and 3.  (In Table 3, there are two factors of $\textsf{P}_n(x)$ which are not listed, both of which have degree $42$.)  Formula (18) implies the Corollary to Theorem 2. \medskip

The class equation $H_{-d}(x)$ for the discriminant $-d \equiv 1$ (mod $8$) may be computed using the resultant
\begin{equation}
cH_{-d}(x)^2=\textrm{Res}_y(b_d(y),y^{16}(1-y^4)x-16(y^8-16y^4+16)^3), \ \ c \in \mathbb{Z}^+,
\end{equation}
at least for small values of the period $n=\textrm{ord}(\tau_d)$ (see the proof of Theorem 6).  This gives a purely algebraic method for calculating $H_{-d}(x)$.  \bigskip

\begin{table}
  \centering 
  \caption{Polynomials $\textsf{P}_n(x), \ 1 \le n \le 5$. \medskip} \label{ }

\begin{tabular}{|c|c|c|c|}
\hline
  &  \\
$n$	&   {$\textsf{P}_n(x)=\prod_{d \in \mathfrak{D}_n}{b_d(x)}$}  \\
  &  \\
\hline
   &  \\
1  &  $x(x+1)(x^2-x+2)$ \ \ ($d=7$) \\
&  \\
2 &   $(x^4-4x^3+5x^2-2x+4)$ \ \ ($d=15$) \\
& \\
3 &  $(x^6+x^5+9x^4-13x^3+18x^2-16x+8)(x^6+7x^5+11x^4-15x^3+16x^2-20x+8)$  \\
& \ \ ($d=23, 31$) \\
& \\
 4  & $(x^8-6x^7+42x^6-60x^5+53x^4-54x^3+24x^2+16)$ \\
 & $\times  (x^8+6x^7+78x^6-84x^5+53x^4-66x^3-12x^2+24x+16)$ \\
 & $\times  (x^8+20x^7+110x^6-100x^5+49x^4-80x^3-40x^2+40x+16)$ \\
 & \ \ \ ($d=39, 55, 63$) \\
 & \\
5 & $(x^{10}-15x^9+74x^8-90x^7+93x^6-187x^5+160x^4-156x^3+168x^2-48x+32)$ \\
 & $\times (x^{10}-31x^9+290x^8-186x^7+5x^6-251x^5-56x^4-60x^3+256x^2+32x+32)$ \\
 & $\times (x^{10}-21x^9+732x^8-290x^7-191x^6-369x^5-502x^4+40x^3+456x^2+144x+32)$ \\
 & $\times (x^{10}+77x^9+1730x^8-366x^7-643x^6-647x^5-1496x^4+120x^3+904x^2+320x+32)$ \\
 & $\times (x^{20}+22x^{19}+1177x^{18}-7012x^{17}+27294x^{16}-72516x^{15}+149882x^{14}-227360x^{13}$ \\
 & $+282149x^{12}-253514x^{11}+152221x^{10}-21772x^9-74372x^8+82952x^7-49328x^6$ \\
 & $-11392x^5+26304x^4-18816x^3+8448x^2+5632x+1024)$ \\
  & \ \ ($d= 47, 79, 103, 127, 119$) \\ 
  & \\
\hline
\end{tabular}

\end{table}

We now prove Theorem 3. \medskip

\noindent {\it Proof of Theorem 3.} \smallskip

The above proof shows that the elements $\pi^4$, where $\pi$ runs through the roots of $b_d(x)$ for which $\pi \cong \pi_d$, are all periodic points of the function $T(z)$ in $\textsf{D}$, and therefore also periodic points of the multivalued function $\hat T(z)$.  We consider the iterated resultants defined in Section 3, but with the polynomial $g(x,y)$ replaced by the polynomial $\tilde g(x,y)=x^2y^2-2(x^2-8x+8)y+x^2$.  This gives us a set of polynomials $\hat R_n(x) \in \mathbb{Z}[x]$, whose roots are the periodic points of $\hat T(z)$ in $\overline{\mathbb{Q}}_2$.  Since the elements $\pi^4$, for $d \in \mathfrak{D}_n$, are all periodic points of $\hat T(z)$, we know that their minimal polynomials divide $\hat R_n(x)$, for any $n \ge 1$.  (Any algebraic conjugate over $\mathbb{Q}$ of a periodic point is also a periodic point.)  Moreover, $\mathbb{Q}(\pi^4)=\Omega_f$, for each such $d$, by \cite{lm}, Proposition 8.1, so that the degree of the minimal polynomial of $\pi^4$ is just $\textrm{deg}(b_d(x))=[\Omega_f:\mathbb{Q}]$.  Noting (19), and that $z=0$ and $z=1$ are fixed points of $\hat T(z)$, this shows that $\hat T(z)$ has at least
$$2+\sum_{k, d}{\textrm{deg}(b_d(x))} \ \ (\textrm{the sum is over} \ k \mid n \ \textrm{and} \ -d \in \mathfrak{D}_k)$$
periodic points in $\overline{\mathbb{Q}}_2$ whose periods divide $n$.  However, similar arguments as in the Lemma of \cite{m2}(pp. 727-728), show that $\textrm{deg}(\hat R_n(x))=2^{n+1}$.  Since
$$2+\sum_{k,d}{\textrm{deg}(b_d(x))}=\textrm{deg}(R_n(x))=2^{n+1},$$
it follows that every root of $\hat R_n(x)$ is $\pi^4$ or $\xi^4$, for some root $\pi$ or $\xi$ of a suitable $b_d(x)$.  This proves Theorem 3.  $\square$. \bigskip

Now we put forward the following theorem and conjecture concerning the discriminants of the polynomials $b_d(x)$. \bigskip

\noindent {\bf Theorem 10.} {\it Let $d$ be any positive integer with $d \equiv 7$ (mod $8$).}
\begin{enumerate}[(a)]
\item {\it If $d> 7$,  the odd prime factors of $\textrm{disc}(b_d(x))$ also divide $\textrm{disc}(H_{-d}(x))$, where $H_{-d}(x)$ is the corresponding class equation.}
\item {\it Any odd prime $p$ which divides $\textrm{disc}(b_d(x))$ satisfies $\left(\frac{-d}{p}\right) \neq 1$.}
\item {\it If $h=h(-d)$, then $2^{3h(h-1)} \mid \textrm{disc}(b_d(x))$.}
\end{enumerate}

\noindent {\it Proof.} From the proof of Theorem 6, the function
$$J(x)=2^8 \frac{(x^2-x+1)^3}{x^2(x-1)^2}$$
satisfies $J(\xi^4)=j(w/2)$, for some ideal basis quotient $w$.  Since $J(x)$ is the $j$-invariant of an elliptic curve in Legendre normal form, we know that $J(x)=J(1-x)$, so
$$J(\pi^4)=J(1-\xi^4)=J(\xi^4)=j(w/2).$$
Conjugating by automorphisms of $\Omega_f/K$ shows that the roots $j_i$ of $H_{-d}(x)$ are given as $j_i=J(\pi_i^4)=J(\xi_i^4)$ for two roots $\pi_i, \xi_i$ of $b_d(x)$ related by $\pi_i^4+\xi_i^4=1$, for $1 \le i \le h(-d)$.  If an odd prime $p$ divides $\textrm{disc}(b_d(x))$, then there is a prime divisor $\mathfrak{p}$ of $p$ in $R_{\Omega_f}$ and either: (i) two roots $\xi_1 \neq \xi_2$ of $b_d(x)$, for which
$$\xi_1 \equiv \xi_2 \ (\textrm{mod} \ \mathfrak{p});$$
(ii) two roots $\xi_1,\pi_2$ of $b_d(x)$ for which
$$\xi_1 \equiv \pi_2 \ (\textrm{mod} \ \mathfrak{p});$$
or (iii) $\xi_1\equiv \pi_1$ (mod $\mathfrak{p})$.  Since $p$ is odd, the quantities $\xi_i^8(\xi_i^4-1)^2=(\xi_i \pi_i)^8\cong 2^8$ are relatively prime to $\mathfrak{p}$.  It follows that for the corresponding roots $j_1 \neq j_2$ of $H_{-d}(x)$, we have
$$j_1=J(\xi_1^4) \equiv J(\xi_2^4)=j_2 \ (\textrm{mod} \ \mathfrak{p})$$
in case (i); and
$$j_1=J(\xi_1^4) \equiv J(\pi_2^4)=j_2 \ (\textrm{mod} \ \mathfrak{p}),$$
in case (ii), implying in either case that $p$ divides the discriminant of $H_{-d}(x)$.  In case (iii) we have from above that $\xi_1, \pi_1 \not \equiv 1$ (mod $\mathfrak{p}$), since $\xi_1^4, \pi_1^4 \not \equiv 1$ (mod $\mathfrak{p}$).  Hence, (8) implies that
$$\pi_1^{\tau^2}=\frac{\xi_1+1}{\xi_1-1} \equiv \frac{\pi_1+1}{\pi_1-1} = \xi_1^{\tau^{-2}} \ (\textrm{mod} \ \mathfrak{p}).$$
If $\tau^2 \neq \tau^{-2}$, this gives
$$\xi_2= \xi_1^{\tau^{-2}} \equiv \pi_1^{\tau^2}=\pi_3 \ (\textrm{mod} \ \mathfrak{p}),$$
where $j_2 \neq j_3$, and we are in case (ii).  If $\tau^2=\tau^{-2}$, then $\tau^4=1$ and the period of the roots of $b_d(x)$ with respect to $\hat F(z)$ is $n=1, 2$, or $4$.  Using (20) and the polynomials in Table 1, we can check the assertion of (a) directly in the cases $d=15, 39, 55, 63$.  This proves (a).  Part (b) follows from the fact that prime factors of $\textrm{disc}(H_{-d}(x))$ satisfy $\left(\frac{-d}{p}\right) \neq 1$.  See \cite{d3}, p. 78.  Part (c) follows easily from the fact that
$$\textrm{disc}(b_d(x))=\prod_{i < j}{(\xi_i-\xi_j)^2} \prod_{i < j}{(\pi_i-\pi_j)^2} \prod_{i,j}{(\xi_i-\pi_j)^2},$$
as follows.  Note that $\wp_2' \mid (\xi_i-\xi_j)$ for $i \neq j$.  Also, Lemma 4 implies that $\xi_i+1 =\frac{\beta_i+2}{2} \cong \wp_2^2$ for all $i$, so that $\wp_2^2 \wp_2' \mid (\xi_i-\xi_j)$ and $\wp_2^3 \wp_2'^3 \cong 2^3 \mid (\xi_i-\xi_j)(\pi_i-\pi_j)$.  Therefore, the first two terms in the above product are divisible by $2^{3h(h-1)}$, as claimed, while the third term is relatively prime to $2$.  
$\square$  \bigskip

\noindent {\bf Conjecture.} \begin{enumerate}[(a)]
\item {\it If $h=h(-d)$, then the exact power of $2$ dividing $\textrm{disc}(b_d(x))$ is $2^{3h(h-1)}$.}
\item {\it If $p$ is an odd prime dividing $\textrm{disc}(b_d(x))$, then $p \le d$.}
\item {\it If $d$ is not prime, the largest prime factor of $\textrm{disc}(b_d(x))$ has the form $q=d-2^k$ for some $k \ge 1$.}
\end{enumerate}

It would be interesting to know if the precise set of primes dividing $\textrm{disc}(b_d(x))$ can be determined, as in Deuring's paper \cite{d3}, or in the conjectures of Yui and Zagier in \cite{yz}.  Also see \cite{gz} and \cite{lv}; the former paper is the starting point for the conjectures in \cite{yz}.

\begin{table}
  \centering 
  \caption{Irreducible factors $b_d(x)$ of $\textsf{P}_6(x)$. \medskip} \label{ }

\begin{tabular}{|c|c|c|c|}
\hline
  &  \\
$d$	&   {$b_d(x)$}  \\
  &  \\
\hline
   &  \\
   87 & $x^{12}+16x^{11}+395x^{10}+398x^9-357x^8-316x^7-155x^6$ \\
& $-1058x^5+1332x^4-704x^3+800x^2-352x+64$ \\
 &   \\
 135 & $x^{12}-36x^{11}+2271x^{10}+1586x^9-1689x^8-1800x^7-2527x^6$ \\
&  $-2310x^5+2664x^4+832x^3+1296x^2-288x+64$ \\
&  \\
175 & $x^{12}-166x^{11}+8027x^{10}+5200x^9-5565x^8-6446x^7-9659x^6$ \\
&  $-6172x^5+6540x^4+5600x^3+2672x^2-32x+64$ \\
&  \\
207  &  $x^{12}-262x^{11}+20035x^{10}+13096x^9-13397x^8-15878x^7-24435x^6$ \\
&  $-14516x^5+14372x^4+15128x^3+5440x^2+416x+64$ \\
&  \\
247 & $x^{12}+184x^{11}+57491x^{10}+39206x^9-36669x^8-44260x^7-70067x^6$ \\
& $-41690x^5+37644x^4+43072x^3+13616x^2+1472x+64$ \\
&  \\
231 & $x^{24}-160x^{23}+39806x^{22}-404188x^{21}+1735295x^{20}-4082916x^{19}$ \\
& $+6591016x^{18}-7995792x^{17}+7025423x^{16}-3646952x^{15}-2986282x^{14}$ \\
&  $+8218276x^{13}-7410127x^{12}+8124428x^{11}-590812x^{10}-4737592x^9$ \\
& $+2208800x^8-5462688x^7+644992x^6+672768x^5+631808x^4$ \\
& $+875008x^3+496640x^2+53248x+4096$ \\
&  \\
255 & $x^{24}+484x^{23}+67682x^{22}-315500x^{21}+1778351x^{20}-3320880x^{19}$ \\
& $+7580476x^{18}-12603888x^{17}+15479855x^{16}-14728444x^{15}+4226978x^{14}$ \\
& $+12258548x^{13}-20944063x^{12}+22569256x^{11}-11161888x^{10}-5859992x^9$ \\
& $+9241280x^8-9494496x^7+2773504x^6+2227200x^5-1364224x^4$ \\
& $+780800x^3+708608x^2+100352x+4096$ \\
 &  \\
\hline
\end{tabular}

\end{table}

\begin{table}
  \centering 
  \caption{Irreducible factors $b_d(x)$ of $\textsf{P}_7(x)$ with $h(-d)=7, 14$. \medskip} \label{ }

\begin{tabular}{|c|c|c|c|}
\hline
  &  \\
$d$	&   {$b_d(x)$}  \\
  &  \\
\hline
  &  \\
71 & $x^{14}-11x^{13}+195x^{12}-127x^{11}+473x^{10}-593x^9+489x^8-1285x^7+1858x^6$ \\
&  $-2880x^5+3320x^4-2656x^3+1792x^2-576x+128$ \\
&   \\
 151 & $x^{14}+49x^{13}+3947x^{12}+5049x^{11}+1257x^{10}-3585x^9-9591x^8-3357x^7-2286x^6$ \\
&  $+756x^5+9648x^4-5760x^3+5152x^2-1280x+128$ \\
&  \\
223 & $x^{14}+327x^{13}+31533x^{12}+49475x^{11}+3971x^{10}-38331x^9-67753x^8-48623x^7-10688x^6$ \\
& $+36240x^5+40216x^4-4912x^3+10848x^2-2304x+128$ \\
&  \\
343$=7^3$ & $x^{14}+553x^{13}+519827x^{12}+864297x^{11}+22281x^{10}-724017x^9-1048551x^8-982269x^7$ \\
& $-51534x^6+722988x^5+532728x^4+113904x^3+33376x^2-3584x+128$ \\
&  \\
463  &  $x^{14}-4317x^{13}+5455509x^{12}+9135083x^{11}+165107x^{10}-7744779x^9-10913545x^8$ \\
&  $-10577543x^7-330800x^6+7742148x^5+5462032x^4+1438400x^3+169824x^2+2880x+128$ \\
 &  \\
   487 & $x^{14} -2219x^{13}+8414699x^{12}+14095377x^{11}+268377x^{10}-11932257x^9-16807911x^8$ \\
& $-16325397x^7-513342x^6+11923632x^5+8397216x^4+2225952x^3+249088x^2+6784x+128$ \\
& \\
\hline
&  \\
287 & $x^{28}+718x^{27}+151595x^{26}+302396x^{25}-1969799x^{24}+13310626x^{23}+49478315x^{22}$ \\  
 & $-92763048x^{21}+15572619x^{20}-55567582x^{19}-49236615x^{18}+258472956x^{17}$ \\  
 & $-25184053x^{16}+97253374x^{15}-160085295x^{14}-143847472x^{13}-97023632x^{12}$ \\
& $+31900208x^{11}+170255840x^{10}-28494624x^9+144550336x^8-146641664x^7$ \\
& $+52719360x^6-60307968x^5+31151104x^4-7360512x^3+3829760x^2-466944x+16384$ \\
&  \\
391 & $x^{28}-910x^{27}+1396079x^{26}-11190416x^{25}+45948277x^{24}-124180050x^{23}+235719087x^{22}$ \\  
 & $-328250004x^{21}+304829895x^{20}-37280970x^{19}-363512763x^{18}+751810392x^{17}$ \\  
 & $-807755041x^{16}+585000802x^{15}+57581533x^{14}-421649716x^{13}+537990116 x^{12}$ \\
& $-439254264x^{11}-53209920x^{10}-5224128x^9-124251648x^8-70235136x^7$ \\
& $+180393984x^6-52199424x^5+110334976x^4+18845696x^3+8744960x^2-401408x+16384$ \\
&  \\
511 & $x^{28}+6614x^{27}+12795083x^{26}-81961412x^{25}+295814809x^{24}-919556958x^{23}+2515624107x^{22}$ \\  
 & $-3835223880x^{21}+2741257515x^{20}-318564558x^{19}-3878860743x^{18}+9526335516x^{17}$ \\  
 & $-6276227797x^{16}+3048095422x^{15}+1197209809x^{14}-7865407120x^{13}+4568895824x^{12}$ \\
& $-2187610536x^{11}+217110912x^{10}+2125718976x^9-1951319616x^8+601344 x^7$ \\
& $-39389184x^6+61917696x^5+688675840x^4+309923840x^3+42622976x^2+1515520x+16384$ \\
 &  \\
\hline
\end{tabular}

\end{table}

\section{Pre-periodic points of $\hat T(z)$ and $\hat F(z)$.}

\noindent {\bf Lemma 11.} {\it The function $\lambda(z)=\frac{\mathfrak{f}_2^8(z)}{\mathfrak{f}^8(z)}$ satisfies the identity}
\begin{equation}
\lambda^2(z) (\lambda(2z)-1)^2=-16\lambda(2z)(\lambda(z)-1), \ \ \textrm{for} \ \Im(z)>0.
\end{equation}

\noindent {\it Proof.} We will show that
\begin{equation}
f(z)=\frac{(\lambda(2z)-1)^2}{\lambda(2z)}=-16\frac{\lambda(z)-1}{\lambda^2(z)}=g(z), \ \ z \in \mathbb{H}.
\end{equation}
We use the fact that $\lambda(z)$ generates the field of modular functions for the subgroup $G$ of $\Gamma=SL_2(\mathbb{Z})$ which is generated by the substitutions $z \rightarrow z+2$ and $z \rightarrow \frac{z}{1-2z}$.  It is clear that $f(z+1)=f(z)$.  Furthermore,
$$\lambda\left(\frac{z}{1-z}\right)=1-\lambda\left(\frac{z-1}{z}\right)=1-\frac{(\lambda(z)-1)}{\lambda(z)}=\frac{1}{\lambda(z)}$$
implies that
$$f\left(\frac{z}{1-2z}\right)=\frac{\left(\frac{1}{\lambda(2z)}-1\right)^2}{\frac{1}{\lambda(2z)}}=f(z).$$
Hence, $f(z)$ is a rational function of $\lambda(z)$.  From \cite{cha}, p. 116, we take the relations
$$\textrm{lim}_{z \rightarrow \infty}\lambda(zi)=0, \ \ \textrm{lim}_{z \rightarrow 0+} \lambda(zi)=1, \ \ \textrm{lim}_{z \rightarrow 0+}\lambda(1+zi) =-\infty.$$
These facts, together with $f(z+1)=f(z)$, imply that $f(z)$ is analytic and nonzero in $\mathbb{H}$ (since $\lambda(z) \neq 0, 1$ in $\mathbb{H}$) and has finite limits at $z=0, z=1$, while $f(z)$ becomes infinite at $z=\infty i$.  Moreover,
\begin{align*}
\lambda(z)&=\frac{2^4 q^{2/3} \prod_{n=1}^\infty{(1+q^{2n})^8}}{q^{-1/3}\prod_{n=1}^\infty{(1+q^{2n-1})^8}}=16q u_1(q)\\
&=16q(1-8q+44q^2-192q^3+718q^4-2400q^5+\cdots), \ \ q=e^{\pi i z},
\end{align*}
where $u_1(q) \in 1+q\mathbb{Z}[[q]]$ and similarly for $u_2(q)$ below.  It follows that the $q$-expansion of $f(z)$ at $\infty i$ is
\begin{align*}
f(z)&=\frac{(-1+16q^2+\cdots)^2}{16q^2+\cdots}=\frac{1}{16}q^{-2}u_2(q)\\
&=\frac{1}{16}q^{-2}-\frac{3}{2}+\frac{69}{4}q^2-128q^4+\frac{5601}{8}q^6-3072q^8+\frac{23003}{2}q^{10}-38400q^{12}+\cdots.
\end{align*}
Therefore,
$$f(z)=\frac{16}{\lambda^2(z)}+b\frac{1}{\lambda(z)}+c$$
is at most a quadratic polynomial in $1/\lambda(z)$.  Then $\textrm{lim}_{z \rightarrow 0+}f(zi)=0$ implies that
$0=16+b+c$.  Finally, using the fact that $\textrm{lim}_{z \rightarrow 0+}\lambda(1+zi) =-\infty$, we have that
$$0=\textrm{lim}_{z \rightarrow 0+} f(zi)=\textrm{lim}_{z \rightarrow 0+} f(1+zi)=c.$$
Hence, $c=0, b=-16$, giving $f(z)=g(z)$, as claimed.  $\square$. \bigskip

Rewriting the identity (21) gives
$$\lambda^2(z) \lambda^2(2z)-2(\lambda^2(z)-8\lambda(z)+8)\lambda(2z)+\lambda^2(z)=0.$$
This shows that $(x,y)=(\lambda(z),\lambda(2z))$ parametrizes the curve
$$\tilde g(x,y)=x^2y^2-2(x^2-8x+8)y+x^2=0,$$
defined by the minimal polynomial $\tilde g(x,y)$ of $y=T(x)$ over $\mathbb{C}(x)$.  Hence, $\lambda(2z)$ is one of the values of $\hat T(\lambda(z))$.  The form of the polynomial $\tilde g(x,y)$ implies that the other root of $\tilde g(\lambda(z),y)=0$ is $y=\frac{1}{\lambda(2z)}=\lambda \left(\frac{2z}{1-2z}\right)$.  Hence, we have for $z \in \mathbb{H}$ that
\begin{equation}
\hat T(\lambda(z)) \in \{\lambda(2z),\frac{1}{\lambda(2z)}\}=\{\lambda(2z),\lambda \left(\frac{2z}{1-2z}\right)\}=\{\lambda(2z),\lambda \left(\frac{2z}{2z+1}\right)\};
\end{equation}
note that $\lambda \left(\frac{2z}{1-2z}\right)=\lambda \left(\frac{2z}{2z+1}\right)$ follows from the fact that $s(z)=\frac{z}{2z+1} \in G$ satisfies $s(\frac{2z}{1-2z})=\frac{2z}{2z+1}$.  Similarly, the inverse function $\hat T^{-1}$ satisfies
\begin{equation}
\hat T^{-1}(\lambda(z)) \in \{\lambda(\frac{z}{2}), \lambda(\frac{z}{2}+1) \}=  \{\lambda(\frac{z}{2}), \frac{\lambda(\frac{z}{2})}{\lambda(\frac{z}{2})-1} \}, \ \ z \in \mathbb{H},
\end{equation}
since $x=\lambda(z),\lambda(z+1)$ are the two roots of $\tilde g(x,\lambda(2z))=0$.  The two images $\lambda(\frac{z}{2})$ and $\frac{\lambda(\frac{z}{2})}{\lambda(\frac{z}{2})-1}$ in (24) coincide exactly when $\lambda(\frac{z}{2})=2$ (since $\lambda(z)$ never takes the value $0$).  The only root of $\tilde g(2,y)=0$ is $y=-1$, and the only root of $\tilde g(x,-1)=0$ is $x=2$, so that (23) and (24) also hold in these cases. \medskip

A pre-periodic point of $\hat T(z)$ is a number $\rho_j$ for which there exist $\rho_i=\hat T^{j-i}(\rho_j)$ satisfying
$$\tilde g(\rho_j,\rho_{j-1})= \cdots = \tilde g (\rho_1,\rho)=0,$$
where $\rho=0,1$, or $\rho=\pi^4$ and $\pi$ is a root of $b_d(x)$ for some $d$.  If $j$ is minimal, the number $\rho_j$ is a pre-periodic point of {\it level $j$}.  Since
$$\tilde g(x,0)=x^2, \ \ \tilde g(x,1)=16(x-1),$$
there are no pre-periodic points corresponding to the fixed points $\rho=0,1$.  Let $p_d(x)$ be the minimal polynomial of $\pi_d^4$, for $d \equiv 7$ (mod $8$).  A pre-periodic point of of $\hat T$ of level $1$ is a number $\rho_1$ for which
$\tilde g(\rho_1,\pi^4)=0$, for some root $\pi^4$ of $p_d(x)$.   Certainly $\rho_1=\pi^{4\tau^{-1}}$ is one solution, by (9) and (17), if $(\pi)=\wp_2$.  By (24), the other solution of this equation is
$$\rho_1=\frac{\pi^{4\tau^{-1}}}{\pi^{4\tau^{-1}}-1}=-\frac{\pi^{4\tau^{-1}}}{\xi^{4\tau^{-1}}}, \ \ \pi^4+\xi^4=1.$$
Similarly, the solutions of $\tilde g(\rho_1,\xi^4)=\tilde g(\rho_1^{\psi},\pi^4)^\psi=0$ are
$$\rho_1=\xi^{4\tau},-\frac{\xi^{4\tau}}{\pi^{4\tau}}=-\left(\frac{\pi^{4\tau^{-1}}}{\xi^{4\tau^{-1}}}\right)^\psi.$$
The roots $\rho_1$ which are not periodic points are roots of the polynomial
$$s_d^{(1)}(x)=(x-1)^{2h(-d)}p_d\left(\frac{x}{x-1}\right).$$
Hence, the $2h(-d)$ roots of $s_d^{(1)}(x)$, all lying in $\Omega_f=\Omega_{2f}$, are the pre-periodic points $\rho_1$ of $\hat T$ at level $1$ corresponding to roots of $p_d(x)$.\medskip

Pre-periodic points of level $2$ corresponding to $p_d(x)$ are numbers $\rho_2$ for which
$$\tilde g(\rho_2,\rho_1)=0, \ \ \rho_1=-\frac{\pi^4}{\xi^4} \ \ \textrm{or} \ \ -\frac{\xi^4}{\pi^4}.$$
I claim that any solution of this equation is an algebraic integer lying in $\Omega_f(i)$.  Such a root is a solution of
\begin{align*}
0=\xi^8 \tilde g(x,-\frac{\pi^4}{\xi^4})&=\pi^8x^2+2(x^2-8x+8)\pi^4\xi^4+\xi^8 x^2\\
&=(\pi^4+\xi^4)^2x^2-16\pi^4 \xi^4 x+16 \pi^4 \xi^4\\
&=x^2-16\pi^4 \xi^4 x+16 \pi^4 \xi^4\\
&=\pi^8 \tilde g(x,-\frac{\xi^4}{\pi^4}).
\end{align*}
The discriminant of the above quadratic is
$$2^8\pi^8\xi^8 -2^6\pi^4\xi^4=2^6 \pi^4 \xi^4(4\pi^4(1-\pi^4)-1)=-2^6 \pi^4 \xi^4 (2\pi^4-1)^2.$$
This proves the claim, and gives us the formula
$$\rho_2=8\pi^4 \xi^4 \pm 4\pi^2 \xi^2 (2\pi^4-1)i=4\pi^2 \xi^2(2\pi^2 \xi^2 \pm (2\pi^4-1) i) \in \Omega_f(i)=\Omega_{4f}$$
for the pre-periodic points at level $2$.  Note that the norm to $\Omega_f$ of the last quantity in this formula is
$$N_{\Omega_f}(2\pi^2 \xi^2 \pm (2\pi^4-1) i)=4\pi^4 \xi^4+4\pi^8-4\pi^4+1=4\pi^4 (\xi^4+\pi^4)-4\pi^4+1=1.$$
Hence, $\rho_2 \cong 2^4$.  Since there are $2h(-d)$ numbers $\rho_1$, and two values of $\rho_1$ give two values of $\rho_2$, by the above displayed equations, there are in all $2h(-d)$ pre-periodic points $\rho_2$ at level $2$.  Also, the numbers $\rho_1$ are all conjugate to each other over $\mathbb{Q}$, and $x^2-16\pi^4 \xi^4 x+16 \pi^4 \xi^4$ is irreducible over $\Omega_f$ ($i \notin \Omega_f$), so the numbers $\rho_2$ are conjugates over $\mathbb{Q}$, as well. \medskip

As above, $s_d(x)=s_d^{(1)}(x)$ denotes the minimal polynomial over $\mathbb{Q}$ of the numbers $\rho_1$ at level $1$, while $s_d^{(2)}(x)$ is the minimal polynomial of the numbers $\rho_2$ at level $2$.  We define the polynomials $s_d^{(r)}(x)$ inductively:
$$s_d^{(r)}(x)=\textrm{Res}_y(\tilde g(x,y),s_d^{(r-1)}(y)), \ \  r \ge 3.$$
Note that for $r=2$, the corresponding resultant in this formula is actually $c(s_d^{(2)}(x))^2$, which is why we start the inductive definition with $r=3$.  The roots of $s_d^{(r)}(x)$ are exactly the pre-periodic points of $\hat T$ at level $r$.  We will prove the following proposition. \bigskip

\noindent {\bf Proposition 12.} {\it The polynomials $s_d^{(r)}(x)$ are irreducible over $\mathbb{Q}$, for $r \ge 1$.  For $r \ge 1$, any root $\rho_r$ of $s_d^{(r)}(x)$ generates the ring class field $K(\rho_r)=\Omega_{2^r f}$ over $K=\mathbb{Q}(\sqrt{-d})$}.  \medskip

\noindent {\it Proof.} We have proved the first assertion for $r=1, 2$.  The second assertion is also clear for $r=1$, since $K(\rho_1)=K(\frac{\pi^4}{\pi^4-1})=K(\pi^4)=\Omega_f=\Omega_{2f}$.  The remaining assertions of the proposition will be proved by induction on $r$.  Let $\pi^4=\lambda(\frac{w}{2})$, where $w=\frac{v+\sqrt{-d}}{2}$ is given in Theorem 3 and $16 \mid (v^2+d)$.  One of the points $\rho_1$ is given by
$$\rho_1=\frac{\pi^4}{\pi^4-1}=\lambda(\frac{w}{2}+1)=\lambda(\frac{w+2}{2}),$$
by (24).  Using (24) repeatedly, we see that
$$\rho_r=\lambda(\frac{w+2}{2^{r}})$$
is a pre-periodic point of $\hat T$ at level $r$. \medskip

Now, the minimal polynomial of the quadratic irrational $\frac{w+2}{2^r}=\frac{v+4+\sqrt{-d}}{2^{r+1}}$ is given by
$$m_r(x)=2^{2r-1}x^2-2^{r-1}(v+4)x+\frac{(v+4)^2+d}{8},$$
with odd constant term, and
$$\textrm{disc} \ m_r(x)=2^{2r-2}(v+4)^2-2^{2r-2}((v+4)^2+d)=-2^{2r-2}d=-(2^{r-1}f)^2 d_K.$$
It follows that $j(\frac{w+2}{2^r})$ generates $\Omega_{2^{r-1}f}$ over $K$.  Furthermore, since $j(z)$ is a rational function of $\lambda(z)$, it is clear that $K(j(\frac{w+2}{2^r})) \subseteq K(\lambda(\frac{w+2}{2^{r}}))$.  \medskip

Assume inductively that $K(\rho_r)=\Omega_{2^r f}$, for some $r \ge 1$.  If we set $s(z)=\frac{(\lambda(z)-1)^2}{\lambda(z)}$ and $t(z)=\frac{\lambda(z)^2}{\lambda(z)-1}$, then we have the identity
$$j(z)=256 \frac{(s(z)+1)^3}{s(z)} = 256 \frac{(t(z)-1)^3}{t(z)}.$$
Setting $s_r=\frac{(\rho_r-1)^2}{\rho_r}$ and $t_r=\frac{\rho_r^2}{\rho_r-1}$ therefore yields
\begin{equation}
j(\frac{w+2}{2^r})=256 \frac{(s_r+1)^3}{s_r}=256 \frac{(t_r-1)^3}{t_r}.
\end{equation}
However, formula (22) gives that
\begin{equation}
s_r=\frac{(\rho_r-1)^2}{\rho_r}=-16\frac{\rho_{r+1}-1}{\rho_{r+1}^2}=\frac{-16}{t_{r+1}}.
\end{equation}
Putting (26) into (25) yields
\begin{equation*}
j(\frac{w+2}{2^{r+2}})=256 \frac{(t_{r+2}-1)^3}{t_{r+2}}=256 \frac{(\frac{-16}{s_{r+1}}-1)^3}{\frac{-16}{s_{r+1}}}=16\frac{(s_{r+1}+16)^3}{s_{r+1}^2}.
\end{equation*}
This shows that $j(\frac{w+2}{2^{r+2}})$ is a rational function of $\rho_{r+1}$, and therefore
\begin{equation}
\Omega_{2^r f}=K(\rho_r)=K\left(j\left(\frac{w+2}{2^{r+1}}\right)\right)  \subset K\left(j\left(\frac{w+2}{2^{r+2}}\right)\right) \subseteq K(\rho_{r+1}).
\end{equation}
However, $[\Omega_{2^{r+1}f}:\Omega_{2^r f}]=2$, for $r \ge 1$, while $[K(\rho_r,\rho_{r+1}):K(\rho_r)] \le 2$.  This forces $K(\rho_{r+1})=K(j(\frac{w+2}{2^{r+2}}))=\Omega_{2^{r+1}f}$ in (27).
 \medskip

Finally, $[\Omega_{2^{r}f}:K]=2^{r-1}h(-d)$ and $\textrm{deg}(s_d^{(r)})=2^{r-1}h(-d)$ (for $r \ge 2$) imply that
$$\textrm{deg}(s_d^{(r+1)})=2 \cdot \textrm{deg}(s_d^{(r)})=2^r h(-d)=[\Omega_{2^{r+1}f}:K]=[K(\rho_{r+1}):K].$$
This shows that $s_d^{(r+1)}(x)$ is irreducible over $K$ (a fortiori over $\mathbb{Q}$), and completes the proof.  $\square$ 
\bigskip

Theorem 3 and Proposition 12 yield the following result, which shows that Conjectures 1 and 2 of \cite{m2} are true for the prime $p=2$ and the function $T(z)$.  (Conjecture 1 was proved for $p=2$ already in \cite{m2}, for a different algebraic function.  The results of this paper give an alternate proof of this conjecture.) \bigskip

\noindent {\bf Theorem 13.} {\it As above, let $p_d(x)$ be the minimal polynomial over $\mathbb{Q}$ of the number $\pi_d^4$, for $d \equiv 7$ (mod $8$).}
\begin{enumerate}[(a)]
\item {\it Any periodic point, other than $0,1$, contained in $\overline{\mathbb{Q}} \subset \overline{\mathbb{Q}}_2$, of the algebraic function $\hat T$ generates a ring class field of odd conductor over an imaginary quadratic field of the form $K=\mathbb{Q}(\sqrt{-d})$, where $-d \equiv 1$ (mod $8$).  Every ring class field of odd conductor over such a field $K$ is generated by an individual periodic point of $\hat T$.  This point can be chosen to be a periodic point of the single valued function $T(z)$ contained in the domain $\textsf{D} \subset \textsf{K}_2$.}
\item {\it Any pre-periodic point of $\hat T$ in $\overline{\mathbb{Q}} \subset \mathbb{C}$ whose level is at least $2$ generates a ring class field of even conductor over one of the fields $K$ in (a). If the level of the pre-periodic point $\rho_r$ is $r \ge 2$, and $\hat T^r(\rho_r)$ is a root of $p_d(x)$, then the conductor of $K(\rho_r)=\Omega_{2^r f}$ over $K=\mathbb{Q}(\sqrt{-d})$ is exactly divisible by $2^r$.  Every ring class field of even conductor over such a field $K$ is generated by an individual pre-periodic point of $\hat T$.}
\end{enumerate}
\bigskip

The minimal polynomials $p_d(x)$ of the periodic points of $\hat T$ are normal polynomials over $\mathbb{Q}$ (see \cite{lm}), while the polynomials $s_d^{(r)}(x)$ are typically not normal over $\mathbb{Q}$, since $\textrm{deg} \ s_d^{(r)}(x) = 2^{r-1}h(-d)=[\Omega_{2^{r}f}:K]$ is only half of the degree of the ring class field $\Omega_{2^r f}$ over $\mathbb{Q}$.  Thus, the ring class fields in question are the normal closures of the fields $\mathbb{Q}(\rho_r)$ they generate over $\mathbb{Q}$, except when $\mathbb{Q}(\rho_r)$ is normal (and therefore abelian) over $\mathbb{Q}$.  In the latter case, $\Omega_{2^r f}=K(\rho_r)=\mathbb{Q}(\sqrt{-d}, \rho_r)$ is abelian over $\mathbb{Q}$, which implies that the discriminant $-2^{2r}d=-4n$, where $n$ is an idoneal number.  See \cite{co1}, pp. 59-62 and \cite{ka}.  From \cite{ka} (see Theorem 3, Corollary 8, and Corollary 23 in that paper), this situation arises only for $2^{2r-2} d=28, 60, 112, 240$; and in these cases $\Omega_{2^r f}=\mathbb{Q}(\rho_r, \sqrt{-d})=\mathbb{Q}(\rho_r, i)$. \medskip

Note that the roots of $s_d^{(r)}(x)$ are invariant under the map $x \rightarrow \frac{x}{x-1}$, by (24). \medskip

The pre-periodic points of $\hat F(z)$ can now be determined in terms of the pre-periodic points of $\hat T(z)$.  Extending the relation (17) to the conjugate value and working in $\overline{\mathbb{Q}}_2$, we have
$$\hat T(z^4) = \hat F(z)^4, \ \ z \in \overline{\mathbb{Q}}_2,$$
meaning that each value of the left-hand side is a value of the right-hand side, and conversely.  If $\alpha_r$ is a pre-periodic point of $\hat F$ at level $r \ge 1$, where $\hat F^r(\alpha_r)=\pi$ is a root of $b_d(x)$, then successive values of $\hat T$ can be chosen so that
$$\hat T^r(\alpha_r^4) = \hat F^r(\alpha_r)^4 = \pi^4$$
is a root of $p_d(x)$; hence, $\alpha_r^4$ is either periodic or pre-periodic of level {\it at most} $r$ with respect to $\hat T$.  On the other hand, if $\rho_r$ is pre-periodic at level $r$ with respect to $\hat T$, and $\rho_r^{1/4}$ is any fourth root of $\rho_r$, then values of $\hat F$ can be chosen so that
$$\hat F^r(\rho_r^{1/4})^4 = \hat T^r(\rho_r) = \pi^4$$
is some root of $p_d(x)$, and therefore
$$\hat F^r(\rho_r^{1/4}) = \varepsilon \pi, \ \ \varepsilon \in \{\pm1, \pm i\}.$$
Now note that $-\pi$ is pre-periodic of level $1$, since $F(-\pi)=F(\pi)$ is a root of $b_d(x)$; and $i\pi$ is pre-periodic of level $2$, if we extend the definition of $F(z)=\frac{-1+\sqrt{1-z^4}}{z^2}$ to the domain $\overline{\textsf{D}}=\{z \in \textsf{K}_2(i): |z|_2 <1\}$ in the quadratic extension $\textsf{K}_2(i)$.  This is because $F(i\pi)=-F(\pi)$.  Hence, $\rho_r^{1/4}$ is either periodic or pre-periodic of level at most $r+2$ with respect to $\hat F$.  (Similar arguments apply in case $T^r(\rho_r)=\xi^4$, for a root $\xi$ which is a unit in $\textsf{K}_2$.)  Though we have worked $2$-adically in this argument, the conclusions are algebraic in nature, and apply equally in $\mathbb{C}$.  This gives the following. \bigskip

\noindent {\bf Theorem 14.} {\it The pre-periodic points of the multivalued function $\hat F(z)$ are the roots of the polynomials
$$\frac{p_d(x^4)}{b_d(x)}= b_d(-x) b_d(ix) b_d(-ix), \ \ d \equiv 7 \ (\textrm{mod} \ 8),$$
together with the roots of the polynomials $s_d^{(r)}(x^4)$, for $r \ge 1$.  The latter roots coincide with the values
$$\lambda(z)^{1/4}=\varepsilon \frac{\mathfrak{f}_2(z)^2}{\mathfrak{f}(z)^2}, \ z=\frac{w+2}{2^r}=\frac{v+4+\sqrt{-d}}{2^{r+1}}, \ r \ge 1, \ \varepsilon \in \{\pm 1, \pm i\}, \ 16 \mid (v^2+d),$$
and their conjugates over $\mathbb{Q}$.  In particular, the periodic and pre-periodic points of $\hat F(z)$ are given by values of modular functions at imaginary quadratic arguments.}
\bigskip

The roots of the polynomials $s_d^{(r)}(x^4)$ generate the same sequence of ring class fields over $K=\mathbb{Q}(\sqrt{-d})$ as do the roots of $s_d^{(r)}(x)$, as we show in the following theorem. \bigskip

\noindent {\bf Theorem 15.} {\it The polynomials $s_d^{(r)}(x^4)$ are irreducible over $\mathbb{Q}$ and $K(\sqrt[4]{\rho_r})=\Omega_{2^{r+2}f}$, for $r \ge 1$.} \medskip

\noindent {\it Proof.} We start from the factorization
$$ \tilde g(x,y^2) = (xy^2-2(x-2)y+x)(xy^2+2(x-2)y+x).$$
Solving for $x$ in this equation gives
$$x=\frac{-4y}{(-y+1)^2} \ \ \textrm{or} \ \ \frac{4y}{(y+1)^2}.$$
Since $\tilde g(\rho_{r+1},\rho_r)=0$, this gives
\begin{equation}
\rho_{r+1}=\frac{4\sqrt{\rho_r}}{(\sqrt{\rho_r}+1)^2}, \ \ r \ge 1,
\end{equation}
for some choice of the square-root.  This implies that $K(\rho_{r+1}) \subseteq K(\sqrt{\rho_r})$, for $r \ge 1$.  Taking the square-root in (28) gives further that
$$\pm \sqrt{\rho_{r+1}}=\frac{2 \sqrt[4]{\rho_r}}{\sqrt{\rho_r}+1}.$$
This implies that
$$\Omega_{2^{r+2}f} = K(\rho_{r+2}) \subseteq K(\sqrt{\rho_{r+1}}) \subseteq K(\sqrt[4]{\rho_r}), \ \ r \ge 1.$$ \smallskip

\noindent Now $[K(\sqrt[4]{\rho_r}):\Omega_{2^r f}]=[K(\sqrt[4]{\rho_r}):K(\rho_r)] \le 4$ and $[\Omega_{2^{r+2}f}:\Omega_{2^r f}]=4$ give that $K(\sqrt[4]{\rho_r})=\Omega_{2^{r+2}f}$.  Since $\textrm{deg}(s_d^{(r)}(x^4))=2^{r+1}h(-d)=[\Omega_{2^{r+2}f}:K]$, this proves that $s_d^{(r)}(x^4)$ is irreducible over $K$, and therefore over $\mathbb{Q}$. $\square$ \bigskip

\noindent {\bf Corollary 16.} {\it (a) The ring class fields $\Omega_{2^r f}$ of even conductor over $K=\mathbb{Q}(\sqrt{-d})$, where $-d \equiv 1$ (mod $8$), are generated by individual pre-periodic points of the algebraic function $\hat F(z)$.  Every pre-periodic point of $\hat F(z)$ whose level is $r \ge 2$ generates a ring class field whose conductor over $K$ is exactly divisible by $2^r$. \smallskip

\noindent (b) All the pre-periodic points $\varepsilon_r$ of $\hat F$ at level $r \ge 1$, for which $\hat F^r(\varepsilon_r)$ is a root of $b_d(x)$, are conjugate over $\mathbb{Q}$.} \medskip 

\noindent {Proof.} We note that the pre-periodic points $\varepsilon_1$ of $\hat F$ at level $1$ are the roots of $b_d(-x)=0$, and the pre-periodic points $\varepsilon_2$ are the roots of $b_d(ix)b_d(-ix)=0$.  In the latter case, the quantities $i\pi$ and $i\xi$ generate $\Omega_{4f}$ over $K$, since $K(\pi)=K(\pi^4)=\Omega_f$ and $\Omega_f(i)=\Omega_{4f}$.  The assertions for $r \ge 3$ will follow directly from the theorem, if we show that the roots $\varepsilon_r=\rho_r^{1/4}$ of $s_d^{(r)}(x^4)$ are pre-periodic of level $r+2$ for the function $\hat F$, for $ r \ge 1$.  We have shown above that the level of the pre-periodic point $\rho_r^{1/4}$ is at most $r+2$.  If it were strictly less than $r+2$, then either $\hat F^{r}(\rho_r^{1/4})= \pi$ or $\xi$ would be a root of $b_d(x)$; or $\hat F^{r}(\rho_r^{1/4})= -\pi$ or $-\xi$.  In the former case, $\hat F^{r-1}(\rho_r^{1/4})=-\pi'$ or $-\xi'$ would be pre-periodic of level $1$. In the latter case, $\hat F^{r-1}(\rho_r^{1/4})=\pm i \pi'$ or $\pm i \xi'$.  In either case, $\rho_r=(\rho_r^{1/4})^4$ would be a pre-periodic point of $\hat T$ of level $\le r-1$.  But this is impossible, since $\rho_r$ is not a root of $s_d^{(i)}(x)$, for $i \le r-1$.  This argument holds for $r \ge 2$, and also for $r=1$, if we set $s_d^{(0)}(x)=p_d(x)$.  This implies all the assertions.  $\square$ \bigskip

These results show that Conjectures 1 and 2 of \cite{m2} (for $p=2$) are true for both functions $\hat T$ and $\hat F$. Since the normal closures of the fields $\mathbb{Q}(\rho_r^{1/4})$ contain the number $i=\sqrt{-1}$, the four cases $d=28, 60, 112, 240$ mentioned above are no longer exceptional in the following theorem.  \bigskip

\noindent {\bf Theorem 17.} {\it The collection of ring class fields over imaginary quadratic fields of the form $K=\mathbb{Q}(\sqrt{-d})$, with $-d \equiv 1$ (mod $8$), coincides with the set of normal closures of the fields generated over $\mathbb{Q}$ by individual periodic and pre-periodic points (different from  $0,-1$) of the algebraic function $\hat F(z)=\frac{-1 \pm \sqrt{1-z^4}}{z^2}$.}

\noindent Dept. of Mathematical Sciences \smallskip

\noindent Indiana University - Purdue University at Indianapolis (IUPUI) \smallskip

\noindent 402 N. Blackford St., Indianapolis, IN 46202 \smallskip

\noindent {\it e-mail}: pmorton@iupui.edu

\end{document}